\documentclass{amsart}

 \theoremstyle{definition}
 
 \theoremstyle{remark}
 
\numberwithin{equation}{section}
\usepackage{amsfonts}
\usepackage{graphicx}
\usepackage{hyperref}

\begin{document}

%\preprint{}

\title[Mathematical computability questions for some classes of differential equations]
{Mathematical computability questions for some classes of\\linear
and non-linear differential equations\\originated from Hilbert's
tenth problem}


\author{Tien D$\tilde {\rm u}$ng Kieu}
\email{kieu@swin.edu.au}
% \thanks{}
\address{Centre for Atom Optics and Ultrafast Spectroscopy,
Swinburne University of Technology, Australia}

\maketitle

\begin{abstract}
Inspired by Quantum Mechanics, we reformulate Hilbert's tenth
problem in the domain of integer arithmetics into problems involving
either a set of infinitely-coupled non-linear differential equations
or a class of linear Schr\"odinger equations with some appropriate
time-dependent Hamiltonians. We then raise the questions whether
these two classes of differential equations are computable or not in
some computation models of computable analysis.  These are
non-trivial and important questions given that: (i) not all
computation models of computable analysis are equivalent, unlike the
case with classical recursion theory; (ii) and not all models
necessarily and inevitably reduce computability of real functions to
discrete computations on Turing machines. However unlikely the
positive answers to our computability questions, their existence
should deserve special attention and be satisfactorily settled since
such positive answers may also have interesting logical consequence
back in the classical recursion theory for the Church-Turing thesis.
\end{abstract}
\begin{center}
Monday, July 4, 2005
\end{center}

\section*{Introduction}
Hilbert's tenth problem~\cite{Davis1982, Matiyasevich1993} is
concerned with the availability of a universal procedure to
determine whether an arbitrarily given Diophantine equation
\begin{eqnarray}
D(x_1,\cdots,x_K) &{=}& 0 \label{1}
\end{eqnarray}
has any positive integer solution or not. After more than 70 years
since its inception, it has finally been shown that the problem is
recursively noncomputable since no such universal procedure which is
also classically recursive can exist.  The problem is equivalent to
the Turing halting problem and intimately links to the concept of
effective computability as defined by the Church-Turing thesis in
classical recursion theory.

Nevertheless, we have established elsewhere, through an inspiration
provided by quantum mechanics~\cite{kieu-contphys, kieu-spie,
kieuFull, kieu-intjtheo, kieu-royal, kieuNew}, some surprising
connections of the above problem in number theory to problems over
the continuous variables. In particular, we have been able to
reformulate Hilbert's tenth problem in terms of a set of infinitely
coupled non-linear differential equations for any given Diophantine
equation~\cite{kieu-royal}. Also, through the framework of quantum
adiabatic computation~\cite{qac}, we have also associated Hilbert's
tenth problem with a class of linear Schr\"odinger equations with
appropriate time-dependent Hamiltonians.  It has been proposed then
that a physical implementation of quantum mechanical processes for
these Schr\"odinger equations could provide the physical means to
solve Hilbert's tenth problem~\cite{kieuNew}.  Here, in this paper
we will only concern ourselves with the differential equations as
mathematical objects, however, and will not appeal to any real
physical processes.  The mathematical objects so derived are
extremely valuable as they provide us a direct link between the
well-established classical recursion theory and the infantile
subject of computable analysis, through a known noncomputable
problem in the former theory.

We especially want to raise in this paper the computability
questions for these differential equations in the domain of
computable analysis, outside and encompassing the classical
recursion theory where Hilbert's tenth problem was originally
formulated.

From the recursive noncomputability of Hilbert's tenth problem, one
might conclude that these differential equations should also be
noncomputable in the wider framework. Such hasty conclusion,
however, is not warranted because of several reasons. Firstly,
however likely the case one might expect, it would need to be
established rigorously as a mathematical truth --because the
unsolvability of Hilbert's tenth problem is only established in the
framework of Turing computability, not necessarily in mathematics in
general.

Secondly, there are many computation models in computable
analysis~\cite{Blum1998, Ko1991, Pour-El1989, Traub1988,
Weihrauch2000} but they are not all equivalent.  And it is known
that computability in one model may not be the same in some other
model, see~\cite{Weihrauch2000}, for example, for a brief discussion
comparison of the various models. This situation is in stark
contrast to the classical recursion theory of functions from
$\mathbb N$ to $\mathbb N$. There, many different formulations have
been given (notably by Kleene, Turing, Post, Herbrand/G\"odel,
Markov) but in the end these all lead to the same notion of
computability with the same class of computable functions. Such an
equivalence has led to the postulation and support of the
Church-Turing thesis in that theory.

The computability definitions for a {\em single} real number in
most, if not all, different computation models of computable
analysis are equivalent.  However, for {\em sequences} of reals the
definitions diverge and are not all equivalent. (The discussions on
sequences of real numbers are necessary because of the necessity of
topological notions in analysis.)  This divergence results in the
dependence of the notion of computability on the different
computation models,\footnote{For example, the so-called real-RAM
approach~\cite{Blum1998, Traub1988} introduces computable real
functions directly and borrows only the concept of {\em control
structure} from Turing computation, without any further referencing
to the latter.  The computability notion in this model, as a result,
is different from that of other real computation models that are
based and built from the Turing computation.} or even on different
choices within a single model.\footnote{One famous example is that
the choice of different norms in a Banach space can lead to opposing
conclusions about the computability of solutions of the same wave
equation in three spatial dimensions with computable initial
functions~\cite{Pour-El1989}.}

In view of such an inequivalence of computability in different
approaches, it is not at all a forgone conclusion that the
noncomputablity of the differential equations mentioned above is
trivially the only, inevitable possibility.

In the next Section we briefly present the observation inspired by
quantum mechanics that leads to the connections between Hilbert's
tenth problem with unbounded self-adjoint operators acting on some
infinite-dimensional Hilbert space. From this we then reformulate
Hilbert's tenth problem in terms of a set of infinitely coupled
non-linear differential equations (eqs.~(\ref{10}, \ref{4}) below).
We also propose a so-called {\em continuation procedure} to
approximate some relevant part of its solution; the computability of
the continuation procedure is left as an unanswered question for the
time being.

We then present a class of linear Schr\"odinger equations
(eq.~(\ref{Schroedinger}) below) with a special class of
time-dependent Hamitonians that are intimately connected to
Hilbert's tenth problem.  We call this the {\em dynamical approach},
because of its use of the Schr\"odinger equations, to distinguish it
from the {\em kinematic approach} above from which we derive the set
of non-linear differential equations.  We conclude the paper with
some remarks and a discussion on the possible implications,
including those on the Church-Turing thesis, if the differential
equations are indeed computable in some computation model of
computable analysis.

\section*{Hilbert's tenth problem and unbounded operators in Hilbert spaces}
Given a Diophantine equation with $K$ unknowns $x$'s as in
eq.~(\ref{1}), we can make a connection,
following~\cite{kieu-contphys, kieu-spie, kieuFull, kieu-intjtheo,
kieu-royal, kieuNew}, with the following self-adjoint operator
acting on some appropriate Fock space (a special type of Hilbert
space)
\begin{eqnarray}
H_P &=& \left(D(a^\dagger_1 a_1,\cdots, a^\dagger_K a_K) \right)^2,
\label{hp}
\end{eqnarray}
where
\begin{eqnarray}
[a_j, a^\dagger_k] = \delta_{jk}, && [a_k, a_j] = 0,
\end{eqnarray}
which are usually termed the creation and annihilation operators,
and most commonly seen in text-book treatment of the quantum simple
harmonic oscillators. The Fock space is built out of the ``vacuum"
$\bigotimes_{j=1}^K|0_j\rangle$ by repeating applications of the
creation operators $a^\dagger_j$.

The operator~(\ref{hp}) has non-negative and discrete eigenvalues
$(D(n_1,\cdots,n_K))^2$, with natural numbers $n_1,\cdots,n_K$.
There is an eigenstate $|E_g\rangle$ corresponding to the smallest
eigenvalue $E_g$.\footnote{Assuming that we could always, if we need
to, eliminate any degeneracy of the eigenvalue by, for example,
modifying the Diophantine equation or by adding some small
perturbation terms in $H_P$.} If the self-adjoint operator is
considered as a Hamiltonian for some dynamical process then these
are respectively the ground state and its energy.

{\em It is then clear that the Diophantine equation~(\ref{1}) has at
least one integer solution if and only if $E_g =
\left(D\left(n_1^{(0)},\cdots, n_K^{(0)}\right) \right)^2 = 0$, for
some $K$-tuple of natural numbers $\left(n_1^{(0)},\cdots,
n_K^{(0)}\right)$.}

To sort out this $E_g$ among the infinitely many eigenvalues is
almost an impossible task.  The strategy we will employ, as inspired
by quantum adiabatic processes, is to tag the state $|E_g\rangle$ by
some other known state $|E_I\rangle$ which is the ground state of
some other self-adjoint operator $H_I$, which can be smoothly
deformed to $H_P$ through some continuous parameter $s\in[0,1]$. To
that end, we consider the interpolating operator
\begin{eqnarray}
{\mathfrak H}(s) &=& H_I + f(s)(H_P - H_I),\nonumber\\
&\equiv& H_I + f(s)W, \label{HFull}
\end{eqnarray}
which has an eigenproblem at each instant $s$,
\begin{eqnarray}
[{\mathfrak H}(s) - E_q(s)]|E_q(s)\rangle = 0, &&q = 0, 1, \cdots
\label{eigen} \label{2}
\end{eqnarray}
with the subscript ordering according to the sizes of the eigenvalues, and $f(s)$
some continuous and monotonically increasing function in $[0,1]$
\begin{eqnarray}
f(0) = 0; && f(1) = 1.
\label{f}
\end{eqnarray}
Clearly, $E_0(0) = E_I$ and $E_0(1) = E_g$.

A suitable $H_I$ is, where  $\alpha$'s $\in {\mathbb C}$,
\begin{eqnarray}
H_I = \sum_{i=1}^K \lambda_i(a^\dagger_i -\alpha_i^*)(a_i
-\alpha_i), \label{HI}
\end{eqnarray}
which we will employ from now on. Here, $E_I = 0$ and $|E_I\rangle =
|\alpha_1\cdots\alpha_K\rangle$ is the Cartesian product of the
coherent states
\begin{eqnarray}
|\alpha_i\rangle &=& {\rm
e}^{-\frac{|\alpha_i|^2}{2}}\sum_{n_i=0}^{\infty}
\frac{\alpha_i^n}{\sqrt{n_i!}}|n_i\rangle. \label{coherent}
\end{eqnarray}
where $|n_i\rangle$ are the eigenstates of $a_i^\dagger a_i$ with
eigenvalues $n_i$.  The $\lambda_i$ can be chosen to be rational or
even irrational numbers such that the first order equation
\begin{eqnarray}
\sum_{i=1}^K \lambda_i p_i &=& 0,
\end{eqnarray}
has no integer solutions in $p_i$.  This condition is to ensure that
all the eigenvalues of $H_I$ are non-degenerate, since the
eigenvalues of $H_I$, which are of the form $\sum_{i=1}^K \lambda_i
n_i$, are then easily seen to be unique for different $K$-tuples of
natural numbers $n$'s.

\section*{The spectral flow -- The ``kinematic" approach}
We now derive the differential equations for the tagging connection
between the instantaneous eigenvalues and eigenvectors at different
instant $s$ in~(\ref{eigen}).

Note firstly that, from the normalisation condition $\langle E_q|E_q\rangle =
1$, we can write
\begin{eqnarray}
\langle E_q|\partial_s|E_q\rangle &=& -i\partial_s\phi_q,
\nonumber
\end{eqnarray}
for some real $\phi_q$.  This can be absorbed away with the redefinition
\begin{eqnarray}
{\rm e}^{i\phi_q(s)}|E_q(s)\rangle &\to& |E_q(s)\rangle,
\end{eqnarray}
upon which
\begin{eqnarray}
\langle E_q|\partial_s|E_q\rangle &=& 0.
\label{comp}
\end{eqnarray}

Differentiating~(\ref{eigen}) with respect to $s$ yields
\begin{eqnarray}
[f'(s)W - \partial_s E_q] |E_q\rangle + [{\mathfrak H} -
E_q]\partial_s|E_q\rangle &=& 0. \label{3}
\end{eqnarray}
We next insert the resolution of unity at each instant $s$,
\begin{eqnarray}
{\bf 1} &=& \sum_{m=0}^\infty |E_m(s)\rangle \langle E_m(s)|,
\nonumber
\end{eqnarray}
just after ${\mathfrak H}$ in~(\ref{3}) to get, by virtue
of~(\ref{comp}),
\begin{eqnarray}
E_q\partial_s|E_q\rangle
&=& [f'(s)W - \partial_sE_q]|E_q\rangle +
\sum_{m\not = q}^\infty E_m \langle E_m|\partial_s|E_q\rangle
|E_m\rangle.
\label{7}
\end{eqnarray}
The inner product of the last equation with $|E_l\rangle$ gives
\begin{eqnarray}
(E_q-E_l)\langle E_l|\partial_s|E_q\rangle &=& f'(s)\langle E_l|W|E_q\rangle
-\partial_sE_q\delta_{ql}.
\label{8}
\end{eqnarray}
Thus, for $q\not=l$ this gives the components of
$\partial_s|E_q\rangle$ in $|E_l\rangle$, provided
$E_q(s)\not=E_l(s)$ at any $s\in (0,1)$, a condition whose proof has
been outlined in~\cite{kieuFull}. Consequently, together
with~(\ref{comp}),
\begin{eqnarray}
\partial_s|E_q\rangle &=& f'(s)\sum_{l\not = q}^\infty \frac{\langle
E_l|W|E_q\rangle}{E_q -E_l} |E_l\rangle.
\label{10}
\end{eqnarray}
Also, putting $q=l$ in~(\ref{8}) we have
\begin{eqnarray}
\partial_s E_q(s) &=& f'(s)\langle E_q(s)| W | E_q(s)\rangle.
\label{4}
\end{eqnarray}
Equations~(\ref{10}) and~(\ref{4}) form the set of infinitely
coupled differential equations providing the tagging linkage we have
been looking for.

In this reformulation, the Diophantine equation~(\ref{1}) has an
integer solution if and only if
\begin{eqnarray}
\lim_{s\to 1} E_0(s) &=& 0, \label{answer}
\end{eqnarray}
from the constructively known eigenvalues and eigenstates of $H_I$
as the initial conditions. The limiting process might be necessary
since $H_P$, i.e. ${\mathfrak H}(1)$, may have a degenerate ground
state eigenvalue in general.
%Formally,
%\begin{eqnarray}
%E_q(s) &=& E_q(0) + \langle E_q(0)|W|E_q(0)\rangle s
%+\nonumber\\
%&& +2\sum_{l\not = q}^\infty {\rm Re}\left\{\int_0^s dt \int_0^t d\tau \frac{\langle
%E_l(\tau)|W|E_q(\tau)\rangle\langle E_q(0)|W|E_l(\tau)\rangle}{E_q(\tau) -E_l(\tau)}
%\right\}+
%\label{e0}\\
%&& +\sum_{l,m\not = q}^\infty \int_0^s dt \int_0^t d\tau \int_0^t d\tau'\frac{\langle
%E_q(\tau)|W|E_m(\tau)\rangle\langle E_l(\tau')|W|E_q(\tau')\rangle
%\langle E_m(\tau)|W|E_l(\tau')\rangle}{\left(E_q(\tau) -E_m(\tau)\right)\left(E_q(\tau') -E_l(\tau')
%\right)}.\nonumber
%\end{eqnarray}

\section*{``Kinematic" continuation procedure?}
The non-linear differential equations of the last Section are
infinitely coupled and may not be solved explicitly or computably in
general. However, we are only interested in the ground state
eigenvalue $E_0(1)$ being zero or not.  And since the influence on
the ground state by states having larger and larger indices
diminishes more and more thanks to the denominators in~(\ref{10}),
this information might be derived in some approximation scheme in
which the number of states involved is truncated to a finite number.
The size of the truncation cannot be universal and must of course
depend on the particular Diophantine equation under consideration.

In the below we speculate on an analytic approximation under the
name of {\em continuation procedure}, and we make no claim about its
computability here but leave it as a challenge for the future.
\begin{itemize}
\item Starting from the initial condition comprising of the
constructively known eigenvalues and eigenvectors of $H_I$ at $s=0$,
the differential equations~(\ref{10}) and~\ref{4}) give us the
series expansions
\begin{eqnarray}
|E_q(\epsilon_1)\rangle &=& |E_q(0)\rangle + \epsilon_1
f'(0)\sum_{l\not = q}^{N_1} \frac{\langle
E_l(0)|W|E_q(0)\rangle}{E_q(0) -E_l(0)} |E_l(0)\rangle
+ \mathfrak R_1,\nonumber\\
E_q(\epsilon_1) &=&  E_q(0) + \epsilon_1f'(0)\langle E_q(0)| W |
E_q(0)\rangle + \mathfrak Q_1.\nonumber \label{continuation}
\end{eqnarray}
\item If the remainders $\mathfrak R_1$ and $\mathfrak Q_1$ above are
computable, we would be able to evaluate the radii of convergence in
$s$, which contain $s=\epsilon_1$, and also the truncation size
$N_1$ which determines the accuracy of the expansions..
\item We then proceed to evaluate new
series approximations similar to the ones above but this time
centred at $s=\epsilon_1$.  The new series have new radii of
convergence and a new truncation of $N_2$ eigenvectors previously
approximated at $s=\epsilon_1$. After this step we have then covered
a finite domain in $s$ away from zero, with some computable degree
of accuracy.
\item We keep reiterating this procedure, if possible, to obtain new
remainders and radii of convergence and thus extend the covered
domain in $s$, until we could evaluate the limit $\lim_{s\to
0}E_0(s)$.
\end{itemize}
% We depict the procedure in Fig.~\ref{fig1}.
This is reminiscent of the procedure of analytic continuation of
functions in complex analysis.
%\begin{figure}
%\begin{center}
%\includegraphics{continuation}
%\caption{\label{fig1}Continuation Procedure.}
%\end{center}
%\end{figure}

Note that at each step of the above procedure we only require some
finite truncation $N_i$ for a given accuracy of the series
approximations.  (That accuracy would also ``recursively" determine
the truncations $N_i$ at all previous steps, $j<i$.)  Note also
that, in general, the condition $x = 0$ for a computable real number
$x$ may not be effectively decidable in some computation model. But
here we have the imposed condition that the eigenvalues at $s =1$
must be integer-valued. This additional condition might help making
the equality condition effectively decidable at $s=1$.  That is, we
would only need to approximate $E_0(1)$ by the procedure above up to
some accuracy, say 0.3, which sufficiently enables us to distinguish
different integers, and we would not require infinite precision. The
imposed condition of integer-valued eigenvalues for $H_P$, by
construction, would provide us the built-in infinite precision at no
extra cost!

\section*{Hilbert's tenth and the Schr\"odinger equation --
The ``dynamical" approach} The decision result for Hilbert's tenth
problem can also be encoded in yet another class of linear
differential equations, apart from the class of nonlinear
equations~(\ref{10}, \ref{4}) above. The linear equation is just the
Schr\"odinger equation which captures the dynamics of a proposed
quantum adiabatic algorithm~\cite{kieu-contphys, kieu-spie,
kieuFull, kieu-intjtheo, kieuNew}.

Let $|\psi(t)\rangle$ be the quantum state at time $t$ (of some
quantum system), its time evolution is given by the Schr\"odinger
equation, for $0<t<\tau$,
\begin{eqnarray}
\partial_t |\psi(t)\rangle &=& -i{\mathfrak H}(t/\tau)|\psi(t)\rangle,
\label{Schroedinger}\\
|\psi(0)\rangle &=& |\alpha_1 \cdots\alpha_K\rangle,\nonumber
\end{eqnarray}
where we have chosen the initial state at time $t=0$ to be the
non-degenerate ground state of $H_I$.

Once again we are only interested in the ground state of $\mathfrak
H (1) = H_P$. The quantum adiabatic theorem~\cite{messiah} asserts
that as $\tau\to\infty$ (that is, when the Hamiltonian $\mathfrak H
(t/\tau)$ in eq.~(\ref{Schroedinger}) varies sufficiently slowly in
the time $t$) the state of the above system at $t=\tau$ would be in
the desired ground state, $|\psi(\tau)\rangle \approx |E_g\rangle$,
to any arbitrary degree of precision!  For a given precision,
various versions of the theorem dictate different conditions on
$\tau$ in terms of some intrinsic properties of the eigenvalues and
eigenfunctions of $\mathfrak H(s)$.  Those conditions thus are
highly dependent on the individual Diophantine equation, and hence
are not suitable to the spirit of a {\em universal} procedure
required by Hilbert's tenth problem.

We have proposed a different and universal criterion for the
identification of the ground state $|E_g\rangle$ from
$|\psi(\tau)\rangle$: {\em the Fock state $\left|n_1^{(0)},\cdots,
n_K^{(0)}\right\rangle$ is the ground state $|E_g\rangle$ if it has
an occupation probability greater than one-half}.  That is,
\begin{eqnarray}
\left|\left\langle \psi(\tau)|n_1^{(0)},\cdots,
n_K^{(0)}\right\rangle\right|^2
> 1/2 &\Rightarrow& |n_1^{(0)},\cdots, n_K^{(0)}\rangle =
|E_g\rangle .\label{criterion}
\end{eqnarray}
It should be emphasised here that such a criterion should be taken
at the present time as a postulate as it has only been proved in
some limited settings~\cite{kieuNew}.  With this criterion, we only
need to repeatedly solve the Schr\"odinger~(\ref{Schroedinger}) for
larger and larger $\tau$ each time until the probability condition
is satisfied so that the ground state can be identified. Such a time
$\tau$ can be shown to exist and be finite by the quantum adiabatic
theorem.  More details of these can be found in~\cite{kieuFull,
kieuNew}.

Note also that while the criterion~(\ref{criterion}) may not the
only one suitable for a physical implementation of the Schr\"odinger
equation, it is the only one that we could yet find suitable for the
mathematical discussion of this paper.

Unlike the case of nonlinear differential equations of a previous
Section, here we could try to make use of a powerful computability
result in a computation model which is known as the First Main
Theorem by Pour-El \& Richards~\cite{Pour-El1989}. Essentially, the
theorem asserts that a {\it bounded} linear operator from a Banach
space to a Banach space which maps a computable sequence of spanning
vectors into another computable sequence will also map {\em any}
computable element into another computable element. For the case at
hand, our Schr\"odinger equation defines a linear operator,
\begin{eqnarray}
U(\tau) &=& {\mathfrak T} {\rm exp}\left\{ -i\tau\int_0^1{\mathfrak
H}(s) ds\right\}, \label{unitary}
\end{eqnarray}
where $\mathfrak T$ is the time-ordering symbol, which maps the
initial state to the final state in the same separable Hilbert
space. Now, our initial state $|\alpha_1 \cdots\alpha_K\rangle$ is
computable by construction.  On the other hand, the linear
operator~(\ref{unitary}) coming out of the Schr\"odinger equation
must be unitary and thus be bounded. Hence, the conditions of the
theorem remained to be checked for fulfillment are: (i) which
mathematical norm should be chosen to enable the identification
criterion~(\ref{criterion}), or some other equivalent criterion, for
the ground state at $t=\tau$; and (ii) whether the image of a
particular computable basis is computable or not with this suitably
chosen norm.

\section*{Concluding remarks}
Inspired by quantum mechanics, we have reformulated the question of
solution existence of a Diophantine equation into the question of
certain properties conceived in an infinitely coupled set of
nonlinear differential equations. In words, we encode the answer of
the former question into the smallest eigenvalue and corresponding
eigenvector of a self-adjoint operator whose integer-valued
eigenvalues are bounded from below. To find these eigen-properties
we next deform the operator continuously to another self-adjoint
operator whose spectrum is known. Once the deformation is also
expressible in the form of a set of nonlinearly coupled differential
equations, we could now start from the constructive knowns as a
handle to study the desired unknowns.

In addition, we also explicitly present a class of linear
Schr\"odinger equations whose solutions at some time $\tau$ from
appropriate initial conditions contain the decision results for any
given Diophantine equation.

These reformulations map a noncomputable problem in the domain of
integer arithmetics into the wider framework of computable analysis.
We have given the names, as explained in the Introduction,
``dynamical" and ``kinematic" respectively for the resulting linear
and non-linear equations. In particular, our set of nonlinear
equations would be an important topic for the largely untouched
subject of nonlinear computable analysis.

For the various reasons given in the Introduction Section, the
questions of computability for these differential equations are
non-trivial and important. Towards some answers for these questions,
we have advocated and speculated on an approach based on the First
Main Theorem by Pour-El \& Richards~\cite{Pour-El1989} for the
Schr\"odinger equations, and some other approximation procedure for
the set of nonlinear differential equations.  However, for now,
these computability questions will have to be left as open problems
and challenges.

If in the (admittedly unlikely) event that there exists computation
model (with suitably chosen norm) in which the differential
equations above are computable and that this model {\em could be
restricted and applied to functions from $\mathbb N$ to $\mathbb N$}
then Hilbert's tenth problem would seem to be solvable in integer
arithmetics (as opposed to be solvable in, for example, some
physical quantum adiabatic computation)! And this would entail a
logical breakdown of the Church-Turing thesis. Because the
restriction of such a computable analysis model to the domain of
integers would provide a new notion of effective computability
different from that of the class of Turing computation.

To be sure that this is {\em not} the case, further investigations
are urgently needed to find out (i) whether the above equations are
noncomputable in {\em all} computation models in computable
analysis; or (ii) whether the model that admits such computability,
if exists, {\em cannot} be restricted to integer arithmetics.

\section*{Acknowledgements}  I am indebted to Peter
Hannaford and Alan Head for discussions and support. I would also
like to thank Marian Pour-El for a discussion on the computability
of quantum mechanics and for her gift of the out-of-print
book~\cite{Pour-El1989}.  This work has been supported by the
Swinburne University Strategic Initiatives.

\bibliography{c:/1data_16Apr05/papers/Computability}

\begin{thebibliography}{10}

\bibitem{Blum1998}
Lenore Blum, Felipe Cucker, Michael Shub, and Steve Smale.
\newblock {\em Complexity and Real Computation}.
\newblock Springer-Verlag, New York, 1998.

\bibitem{Davis1982}
Martin Davis.
\newblock {\em Computability and Unsolvability}.
\newblock Dover, New York, 1982.

\bibitem{qac}
E.~Farhi, J.~Goldstone, S.~Gutmann, and M.~Sipser.
\newblock Quantum computation by adiabatic evolution.
\newblock {\tt ArXiv:quant-ph/0001106}, 2000.

\bibitem{kieu-contphys}
T.D. Kieu.
\newblock Computing the non-computable.
\newblock {\em Contemporary Physics}, 44:51--77, 2003.

\bibitem{kieu-spie}
T.D. Kieu.
\newblock Numerical simulations of a quantum algorithm for {H}ilbert's tenth
  problem.
\newblock In Eric Donkor, Andrew~R. Pirich, and Howard~E. Brandt, editors, {\em
  Proceedings of SPIE Vol. 5105 {\it Quantum Information and Computation}},
  pages 89--95. SPIE, Bellingham, WA, 2003.

\bibitem{kieuFull}
T.D. Kieu.
\newblock Quantum adiabatic algorithm for {H}ilbert's tenth problem: I. {T}he
  algorithm.
\newblock \texttt{ArXiv:quant-ph/0310052}, 2003.

\bibitem{kieu-intjtheo}
T.D. Kieu.
\newblock Quantum algorithms for {H}ilbert's tenth problem.
\newblock {\em Int. J. Theor. Phys.}, 42:1451--1468, 2003.

\bibitem{kieu-royal}
T.D. Kieu.
\newblock A reformulation of {H}ilbert's tenth problem through quantum
  mechanics.
\newblock {\em Proc. Roy. Soc.}, A 460:1535--1545, 2004.

\bibitem{kieuNew}
T.D. Kieu.
\newblock Hypercomputability of quantum adiabatic processes: Fact versus
  prejudices.
\newblock \texttt{ArXiv:quant-ph/0504101}, an invited paper for and to appear
  in a special issue of the Journal of Applied Mathematics and Computation,
  2005.

\bibitem{Ko1991}
Ker-I Ko.
\newblock {\em Complexity Theory of Real Functions}.
\newblock Birkh\"aser, Boston, 1991.

\bibitem{Matiyasevich1993}
Yuri~V. Matiyasevich.
\newblock {\em Hilbert's Tenth Problem}.
\newblock MIT Press, Cambridge, 1993.

\bibitem{messiah}
A.~Messiah.
\newblock {\em Quantum Mechanics}.
\newblock Dover, New York, 1999.

\bibitem{Pour-El1989}
Marian~B. Pour-El and J.~Ian Richards.
\newblock {\em Computability in Analysis and Physics}.
\newblock Springer-Verlag, Berlin Heidelberg, 1989.

\bibitem{Traub1988}
Joseph~F. Traub, G.W. Wasilkowski, and H.~Wozniakowski.
\newblock {\em Information-Based Complexity}.
\newblock Acaemic Press, New York, 1988.

\bibitem{Weihrauch2000}
Klaus Weihrauch.
\newblock {\em Computable Analysis}.
\newblock Springer-Verlag, Berlin Heidelberg, 2000.

\end{thebibliography}
\bibliographystyle{plain}
%Included for Gather Purpose only:
%input "c:/1data_16Apr05/papers/Computability.bib"
%input "adiabatic.bib"
%input "ait.bib"

\end{document}